\documentclass[11pt,a4paper, reqno]{amsart}
\usepackage[T1]{fontenc}
\usepackage[utf8]{inputenc}
\usepackage{lmodern}
\usepackage{amsmath,amssymb,amsthm,mathtools}
\usepackage[margin=1in]{geometry}
\usepackage{microtype}
\usepackage{enumitem}
\usepackage{tikz-cd}
\usepackage{cite}
\usepackage{hyperref}
\hypersetup{colorlinks=true,citecolor=blue,linkcolor=blue,urlcolor=blue}

\newtheorem{theorem}{Theorem}[section]
\newtheorem{proposition}[theorem]{Proposition}
\newtheorem{lemma}[theorem]{Lemma}
\newtheorem{corollary}[theorem]{Corollary}
\theoremstyle{remark}
\newtheorem{remark}[theorem]{Remark}

\newcommand{\I}{[0,1]}
\newcommand{\R}{\mathbb{R}}

\newcommand{\Ci}{C(I)}
\newcommand{\Sp}{S(C^1(I))^+}
\newcommand{\K}{\mathcal K}
\newcommand{\one}{\mathbf{1}}
\newcommand{\dd}{\mathrm{d}}
\newcommand{\norm}[1]{\left\lVert #1\right\rVert}
\newcommand{\norminf}[1]{\left\lVert #1\right\rVert_\infty}

\newcommand{\id}{\operatorname{id}}

\title[Anchored positive-sphere isometries of $C^1$]{The positive-sphere variant of Tingley's problem for $C^1([0,1])$}
\author{Min-Ruei Lin}
\address{Department of Applied Mathematics, National Sun Yat-sen University, Kaohsiung 80424, Taiwan; and Graduate School of Science and Technology, Niigata University, Niigata 950--2181, Japan}
\email{m082030021@student.nsysu.edu.tw}

\begin{document}
\begin{abstract}
Let the real space $C^1([0,1])$ be endowed with the norm $\|f\|=\|f\|_\infty+\|f'\|_\infty$.
For fixed $x_0\in[0,1]$, let $S(C^1([0,1]))^+_{x_0}$ be the part of its unit sphere consisting of nondecreasing functions satisfying $f(x_0)\ge 0$.
For $x_1,x_2\in[0,1]$, we prove that a mapping $\Delta:S(C^1([0,1]))^+_{x_1}\to S(C^1([0,1]))^+_{x_2}$ is a surjective isometry if and only if $x_1=x_2$ and $\Delta$ is the identity mapping.
In this case, the unique real-linear extension of $\Delta$ to $C^1([0,1])$ is the identity operator.
\end{abstract}
\maketitle

\section{Introduction}
Tingley's problem asks whether a surjective isometry between the unit spheres of real Banach spaces extends to a surjective $\R$-linear isometry of the ambient spaces \cite{Tingley}.
For a survey in the setting of operator algebras, see Peralta \cite{Peralta}.
In the positive-sphere variant, one replaces each unit sphere by its intersection with a generating positive cone.
For compact Hausdorff spaces, Leung, Ng and Wong \cite{LeungNgWong} showed that surjective isometries of the positive unit spheres of continuous functions with the supremum norm are induced by homeomorphisms.

For the positive sphere of $C^1[0,1]$, Enami, Hirota, Koshimizu and Lin \cite[Theorem~1.1]{EHKL} obtain a complex-linear extension for the norms $\|(|f(0)|,\norminf{f'})\|_{\ell^p}$, $1\leq p\leq\infty$.
For $p=1$ their function-value term is $|f(0)|$.
In contrast, we use the sum norm involving $\norminf{f}$.
Even for the cone of nonnegative nondecreasing functions, this introduces an integral term in derivative coordinates, both in the sphere constraint and in the transported metric.
A general affirmative solution to the classical Tingley's problem has recently been proposed in the preprint \cite{OpenAITingley}.
Such a result does not directly settle its positive-sphere variant, since a positive-sphere isometry need not a priori extend to an isometry of the full unit sphere.

Put $I=\I$, and let $C^1(I)$ simply denote $C^1(\I,\R)$.
Equip $C^1(I)$ with the norm
\[
\norm{f}:=\norminf{f}+\norminf{f'}.
\]
For a fixed anchor $x_0\in I$, define
\[
\K_{x_0}:=\{f\in C^1(I):f(x_0)\geq0,\ f'\geq0\},\qquad
\Sp_{x_0}:=\{f\in\K_{x_0}:\norm{f}=1\}.
\]
All function inequalities are pointwise, and $\one$ denotes the constant function identically equal to $1$.
It is readily verified that $\K_{x_0}$ is a closed, pointed, convex cone.
To see that $\K_{x_0}$ is generating, take any $f\in C^1(I)$ and define
\[
u(x)=\max\{0,-f(x_0)\}+\norminf{f'}(x-x_0).
\]
Then $u(x_0)\geq0$, $(f+u)(x_0)\geq0$, $u'\geq0$ and $(f+u)'\geq0$.
Hence $u,f+u\in\K_{x_0}$, and the decomposition $f=(f+u)-u$ shows that $C^1(I)\subseteq\K_{x_0}-\K_{x_0}$.
Converse is immediate, so
\[
C^1(I)=\K_{x_0}-\K_{x_0}.
\]
Notice that, for $x_0>0$, there exists $h\in\K_{x_0}$ such that $h(t)<0$ for some $0\le t<x_0$.
Furthermore, monotonicity gives the normalization formula
\begin{equation}\label{eq:normalization-general}
\norm{f}=\max\{-f(0),f(1)\}+\norminf{f'}\qquad(f\in\K_{x_0}).
\end{equation}

Our main result is the following.
\begin{theorem}\label{thm:positive-rigidity}
Let $x_1,x_2\in I$.
If $\Delta:\Sp_{x_1}\to\Sp_{x_2}$ is a surjective isometry, then $x_1=x_2=:x_0$ and $\Delta=\id_{\Sp_{x_0}}$.
Moreover, $\Delta$ admits a unique $\R$-linear extension to $C^1(I)$, namely $\id_{C^1(I)}$.
\end{theorem}

The change from the endpoint anchor to an arbitrary anchor is not merely a change of coordinates.
The two terms in the maximum in \eqref{eq:normalization-general} give rise to two branches of the sphere, and the derivative map need not be injective on the entire sphere.
Our proof first fixes $\one$ and identifies the branch on which the terminal value strictly dominates the negative initial value.
Within that branch, a metric characterization of constant shifts by unique midpoints permits the recovery of the supremum-norm metric.
We then use Mankiewicz's theorem and the real Banach--Stone theorem, and recover the other branch from metric intervals.

\section{Metric identification of a branch and derivative coordinates}
Set $\Ci=C(I,\R)$ and define $Df=f'$ on $C^1(I)$.
Fix $x_0\in I$ and a surjective isometry $\Delta:\Sp_{x_0}\to \Sp_{x_0}$.
For $f\in \Sp_{x_0}$, since $f$ is nondecreasing and $f(x_0)\geq0$, we have 
\begin{equation}\label{eq:distance-one}
\norm{\one-f}=1-f(0)+\norminf{f'}.
\end{equation}

\begin{lemma}\label{lem:constant-fixed}
The isometry $\Delta$ fixes $\one$.
More precisely, $\one$ is the unique point $f\in \Sp_0$ with $\sup_{g\in \Sp_0}\norm{f-g}=2$.
If $x_0>0$, then $\one$ is the unique point of $\Sp_{x_0}$ having at least two distinct points at distance $2$ from it.
\end{lemma}
\begin{proof}
Suppose first that $x_0=0$.
In this case all functions in $\Sp_0$ are nonnegative, and \eqref{eq:normalization-general} gives $f(1)=1-\norminf{f'}$.
For $f,g\in \Sp_0$, nonnegativity of the functions and their derivatives gives
\[
\norm{f-g}\leq \max\{1-\norminf{f'},1-\norminf{g'}\}+\max\{\norminf{f'},\norminf{g'}\}=1+|\norminf{f'}-\norminf{g'}|.
\]
We have $\norminf{f'}<1$ for all $f\in \Sp_0$, since otherwise $f(1)=0$ would force $f=0$.
If $f\ne\one$, then also $\norminf{f'}>0$, so
\[
\sup_{g\in \Sp_0}\norm{f-g}\le 1+\max\{\norminf{f'},1-\norminf{f'}\}<2.
\]
On the other hand, for $j\ge1$ set
\[
g_j(x)=\frac{\int_0^x t^j\,\dd t}{1+\int_0^1 t^j\,\dd t}.
\]
Then $g_j\in \Sp_0$ and, by \eqref{eq:distance-one}, $\norm{\one-g_j}=1+(1+1/(j+1))^{-1}\to2$.
This proves the first characterization.

Now suppose $x_0>0$.
For arbitrary $f,g\in \Sp_{x_0}$, the nonnegativity of their derivatives implies
\begin{align*}
\norm{f-g}&\leq \norminf{f}+\norminf{g}+\max\{\norminf{f'},\norminf{g'}\}
=2-\min\{\norminf{f'},\norminf{g'}\}.
\end{align*}
Since $\one$ is the only constant member of $\Sp_{x_0}$, any pair of points at distance $2$ must contain $\one$.
There are at least two such pairs.
To see this, for any nonzero $h\in\Ci$ with $h\geq0$ and $h=0$ on $[x_0,1]$, define
\[
f_h(x)=\frac{\int_{x_0}^x h(t)\,\dd t}{\norminf{h}+\int_0^{x_0} h(t)\,\dd t}.
\]
Then $f_h\in \Sp_{x_0}$, $f_h(1)=0$ and $-f_h(0)=\norminf{f_h}$.
Equation~\eqref{eq:distance-one} gives $\norm{\one-f_h}=2$.
Choosing two nonproportional functions $h$ supported in $(0,x_0)$ gives distinct points $f_h$.
This proves the second characterization.
Both characterizations are preserved by surjective isometries.
\end{proof}
In particular, $\Sp_{0}$ has diameter $2$, although its diameter is not attained by any pair of points.

Define the two subsets
\begin{equation}\label{eq:branches}
U_{x_0}:=\{f\in \Sp_{x_0}:f(1)>-f(0)\},\qquad E_{x_0}:=\Sp_{x_0}\setminus U_{x_0}.
\end{equation}
We refer to $U_{x_0}$ as the strict right branch; the set $E_{x_0}$ includes the common boundary where $f(1)=-f(0)$.

\begin{lemma}\label{lem:right-branch}
One has $U_{x_0}=\{f\in \Sp_{x_0}:\norm{\one-f}<2\}$.
In particular, $\Delta(U_{x_0})=U_{x_0}$ and $\Delta(E_{x_0})=E_{x_0}$.
\end{lemma}
\begin{proof}
By \eqref{eq:normalization-general} and \eqref{eq:distance-one}, $\norm{\one-f}=2$ if and only if $-f(0)=\norminf{f}$, or equivalently $-f(0)\geq f(1)$.
The assertion follows from Lemma~\ref{lem:constant-fixed}.
\end{proof}

For nonnegative $h\in\Ci$ we have $\norminf{h}=\max_I h$.
Introduce the convex set
\begin{equation}\label{eq:derivative-domain}
K_{x_0}:=\left\{h\in\Ci:\begin{array}{l}
h\geq0,\quad \norminf{h}+\int_{x_0}^1h(t)\,\dd t\leq1,\quad
2\norminf{h}+\int_0^1h(t)\,\dd t<2
\end{array}\right\}
\end{equation}
and its closed counterpart
\begin{equation}\label{eq:K-general}
\overline{K_{x_0}}:=\left\{h\in\Ci:\begin{array}{l}
h\geq0,\quad \norminf{h}+\int_{x_0}^1h(t)\,\dd t\leq1,\quad
2\norminf{h}+\int_0^1h(t)\,\dd t\leq2
\end{array}\right\}.
\end{equation}

\begin{proposition}\label{prop:right-derivative-bijection}
The derivative map $D:f\mapsto f'$ restricts to a bijection from $U_{x_0}$ onto $K_{x_0}$ with inverse
\begin{equation}\label{eq:Phi-general}
\Phi(h)(x)=1-\norminf{h}-\int_x^1 h(t)\,\dd t.
\end{equation}
Moreover, $\overline{K_{x_0}}$ is a closed convex body in $\Ci$ with nonempty interior, and $K_{x_0}$ is dense in $\overline{K_{x_0}}$ for the supremum norm.
\end{proposition}
\begin{proof}
If $f\in U_{x_0}$ and $h=f'$, then $f(1)=1-\norminf{h}$ by \eqref{eq:normalization-general}, whence $f=\Phi(h)$.
The condition $f(x_0)\geq0$ is equivalent to $\norminf{h}+\int_{x_0}^1h(t)\,\dd t\leq1$, and $f(1)>-f(0)$ is equivalent to $2\norminf{h}+\int_0^1h(t)\,\dd t<2$.
Conversely, for $h\in K_{x_0}$, these same relations give $\Phi(h)(x_0)\geq0$, $\Phi(h)'=h\geq0$ and $\Phi(h)(1)>-\Phi(h)(0)$.
Thus $\norm{\Phi(h)}=1$ and $\Phi(h)\in U_{x_0}$.

The convexity and closedness of $\overline{K_{x_0}}$ follow from continuity and convexity of the norm and linearity of the integrals.
It is bounded, since $\norminf{h}\leq1$ on $\overline{K_{x_0}}$.
The open supremum-norm ball with center $\one/4$ and radius $1/8$ lies in $K_{x_0}$: every member is strictly positive and has norm less than $3/8$, hence satisfies both constraints strictly.
Finally, $(1-1/j)h\in K_{x_0}$ for every $h\in \overline{K_{x_0}}$ and $j\geq2$, which proves the density assertion.
\end{proof}

\begin{remark}\label{rem:two-branches}
For $f\in E_{x_0}$ with $h=f'$, normalization instead gives
\begin{equation}\label{eq:Phi-minus}
f(x)=\norminf{h}-1+\int_0^x h(t)\,\dd t,\qquad 2\norminf{h}+\int_0^1h(t)\,\dd t\leq2.
\end{equation}
The derivative map on $\Sp_{x_0}$ is not injective whenever $x_0>0$.
In fact, choose a nonzero nonnegative continuous $q$ supported in $(0,x_0)$.
The function
\[
f(x)=\frac{\int_{x_0}^x q(t)\,\dd t}{\norminf{q}+\int_0^1q(t)\,\dd t}
\]
lies in $E_{x_0}$ with $f(1)=0$ and $-f(0)>0$.
Its derivative
\[
h=\frac{q}{\norminf{q}+\int_0^1q(t)\,\dd t}
\]
satisfies $2\norminf{h}+\int_0^1h(t)\,\dd t<2$, and $\Phi(h)=f+(2-2\norminf{h}-\int_0^1h(t)\,\dd t)\one$ belongs to $U_{x_0}$ and has the same derivative.
By contrast, for $x_0=0$ we have $E_0=\varnothing$, since $f(0),f(1)\geq0$ and $f(1)=0$ is impossible for a unit vector.
\end{remark}

Transport the metric of $U_{x_0}$ to $K_{x_0}$ via $\Phi$.
With the linear operator $J$ given by $(Jv)(x):=\int_x^1v(t)\,\dd t$, the metric is
\begin{equation}\label{eq:dK}
\begin{aligned}
d_K(h,k)&:=\norm{\Phi(h)-\Phi(k)}\\
&=\norminf{h-k}+\norminf{(\norminf{h}-\norminf{k})\one+J(h-k)}.
\end{aligned}
\end{equation}
We call an isometry with respect to the transported metric $d_K$ a $d_K$-isometry.
By Lemma~\ref{lem:right-branch}, the map
\begin{equation}\label{eq:T}
T:=D\circ\Delta\circ\Phi:K_{x_0}\longrightarrow K_{x_0}
\end{equation}
is a bijective $d_K$-isometry.
Since $D|_{U_{x_0}}$ and $\Phi$ are mutually inverse, the following diagram commutes.
\[
\begin{tikzcd}[row sep=large, column sep=huge]
U_{x_0}
    \arrow[r, "\Delta"]
    \arrow[d, "D"', shift right=1]
&
U_{x_0}
    \arrow[d, "D"', shift right=1]
\\
K_{x_0}
    \arrow[r, "T"']
    \arrow[u, "\Phi"', shift right=1]
&
K_{x_0}
    \arrow[u, "\Phi"', shift right=1]
\end{tikzcd}
\]
In addition, $T0=0$ because $\Phi(0)=\one$.
For all $h,k\in K_{x_0}$,
\begin{equation}\label{eq:metric-comparison}
\norminf{h-k}\leq d_K(h,k)\leq3\norminf{h-k},
\end{equation}
since $|\norminf{h}-\norminf{k}|\leq\norminf{h-k}$ and $\norminf{J(h-k)}\leq\norminf{h-k}$.
In particular $T$ and $T^{-1}$ are continuous in the supremum norm.
Taking $k=0$ in \eqref{eq:dK} and using $h\geq0$ yields
\begin{equation}\label{eq:R-invariant}
\begin{aligned}
d_K(0,h)&=2\norminf{h}+\int_0^1h(t)\,\dd t,\\
2\norminf{Th}+\int_0^1(Th)(t)\,\dd t&=2\norminf{h}+\int_0^1h(t)\,\dd t.
\end{aligned}
\end{equation}

\section{Unique metric midpoints and admissible constant shifts}\label{sec:midpoints}
For $f,g\in \Sp_{x_0}$, a \emph{metric midpoint} means a point $z\in \Sp_{x_0}$ satisfying $\norm{f-z}=\norm{g-z}=\norm{f-g}/2$.
No convexity of the sphere is assumed in this definition.
The following lemma provides an intrinsic characterization of constant shifts in derivative coordinates.

\begin{lemma}\label{lem:unique-midpoint}
Let $f,g\in U_{x_0}$ with $f\ne g$.
The following are equivalent:
\begin{enumerate}
\item There is $c>0$ such that $g'=f'+c\one$.
\item The equality
\begin{equation}\label{eq:metric-interval}
\norm{\one-g}=\norm{\one-f}+\norm{f-g}
\end{equation}
holds, and $f,g$ have a unique metric midpoint in $\Sp_{x_0}$.
\end{enumerate}
\end{lemma}
\begin{proof}
Suppose first that $g'=f'+c\one$ with $c>0$.
Formula~\eqref{eq:Phi-general} gives
\begin{equation}\label{eq:function-shift}
g(x)-f(x)=-2c+cx.
\end{equation}
Thus $\norm{f-g}=3c$, and \eqref{eq:distance-one} implies \eqref{eq:metric-interval}.
Since $K_{x_0}$ is convex and $\norminf{f'+s\one}=\norminf{f'}+s$ for $s\geq0$, the arithmetic midpoint $w=(f+g)/2=\Phi(f'+(c/2)\one)$ belongs to $U_{x_0}$ and is a metric midpoint.

Let $z$ be any metric midpoint in $\Sp_{x_0}$.
Equality in the triangle inequality for the sum norm forces equality in both component triangle inequalities.
In particular, $\norminf{z'-f'}+\norminf{g'-z'}=c$.
At each $x\in I$,
\[
c=(z'-f')(x)+(g'-z')(x)\leq\norminf{z'-f'}+\norminf{g'-z'}=c.
\]
Hence $z'=f'+\norminf{z'-f'}\one$.
By \eqref{eq:metric-interval} and Lemma~\ref{lem:right-branch},
\[
\norm{\one-z}\leq\norm{\one-f}+\tfrac12\norm{f-g}
=\tfrac12(\norm{\one-f}+\norm{\one-g})<2,
\]
so $z\in U_{x_0}$.
The derivative determines $z$ uniquely there, and \eqref{eq:Phi-general} gives $(z-f)(x)=-2\norminf{z'-f'}+x\norminf{z'-f'}$.
Consequently $3\norminf{z'-f'}=\norm{z-f}=3c/2$ and $z=w$.

Conversely, suppose that \eqref{eq:metric-interval} holds and the metric midpoint is unique.
Write $h=f'$ and $k=g'$.
As $D$ is injective on $U_{x_0}$, we have $\norminf{k-h}>0$.
Equality in the derivative component of \eqref{eq:metric-interval} yields $\norminf{k}=\norminf{h}+\norminf{k-h}$.
Choose $t_0\in I$ with $k(t_0)=\norminf{k}$.
Then
\[
h(t_0)\geq k(t_0)-\norminf{k-h}=\norminf{h},
\]
so $h(t_0)=\norminf{h}$ and $(k-h)(t_0)=\norminf{k-h}$.
In particular,
\[
\norminf{(h+k)/2}=(\norminf{h}+\norminf{k})/2.
\]
Thus $w=(f+g)/2$ is in $U_{x_0}$: it has nonnegative derivative, satisfies the anchor condition and strict right dominance by convexity of these inequalities, and its terminal value plus its derivative norm is one.
It is a metric midpoint.

Suppose that $k-h$ is not identically $\norminf{k-h}$.
By continuity and the intermediate value theorem, the set on which $0<k-h<\norminf{k-h}$ is relatively open and nonempty.
We may therefore choose a nonempty open interval $V$ whose closure is contained in this set and in $(0,1)\setminus\{x_0,t_0\}$.
For $x\in V$ we have
\begin{equation}\label{eq:perturbation-slack}
0<w'(x)<\norminf{w'},\quad |k(x)-h(x)|<\norminf{k-h},\quad
|g(x)-f(x)|<\norminf{g-f}.
\end{equation}
The first strict upper bound follows because equality would require $h(x)=\norminf{h}$ and $k(x)=\norminf{k}$, giving $(k-h)(x)=\norminf{k-h}$.
The last strict inequality follows from $(g-f)'(x)>0$: a differentiable function cannot attain either a global maximum or a global minimum at such an interior point.

Choose a nonzero $\psi\in C_c^1(V)$, extended by zero to $I$, and let $z=w+\varepsilon\psi$.
All strict inequalities in \eqref{eq:perturbation-slack} have uniform positive margins on the compact support of $\psi$.
For sufficiently small $\varepsilon>0$, we therefore have $0\leq z'\leq \norminf{w'}$ everywhere.
Since $z'(t_0)=w'(t_0)=\norminf{w'}$, it follows that $\norminf{z'}=\norminf{w'}$.
The endpoint and anchor values are unchanged; monotonicity and \eqref{eq:normalization-general} now imply $z\in U_{x_0}$.
Also, on the support of $\psi$ the inequalities
\begin{align*}
|z-f|,\ |z-g|&\leq \tfrac12\norminf{f-g},\\
|z'-f'|,\ |z'-g'|&\leq \norminf{k-h}/2
\end{align*}
hold for sufficiently small $\varepsilon$, by the last two margins in \eqref{eq:perturbation-slack}.
Outside the support they hold because $z=w$.
Hence $\norm{z-f},\norm{z-g}\leq\norm{f-g}/2$, and the triangle inequality forces equality in both.
This produces a metric midpoint $z\ne w$, a contradiction.
Therefore $k-h=\norminf{k-h}\one$, as required.
\end{proof}

\begin{proposition}\label{prop:constant-shifts}
If $h\in K_{x_0}$ and $c\geq0$ satisfy $h+c\one\in K_{x_0}$, then
\begin{equation}\label{eq:shift-equivariance}
T(h+c\one)=Th+c\one.
\end{equation}
The same assertion holds for $T^{-1}$.
\end{proposition}
\begin{proof}
The case $c=0$ is trivial.
For $c>0$ apply Lemma~\ref{lem:unique-midpoint} to $f=\Phi(h)$ and $g=\Phi(h+c\one)$.
A surjective isometry preserves metric intervals and bijectively maps the sets of metric midpoints, and $\Delta$ fixes $\one$ and maps $U_{x_0}$ onto itself.
Applying the lemma to $\Delta f,\Delta g$ yields $T(h+c\one)=Th+b\one$ for some $b>0$.
Formula~\eqref{eq:function-shift} gives $d_K(h,h+c\one)=3c$ and $d_K(Th,Th+b\one)=3b$.
Hence $b=c$.
The proof for $T^{-1}$ is identical, with $\Delta^{-1}$ in place of $\Delta$.
\end{proof}

\section{Recovery of the supremum-norm metric}
Although $T$ is a $d_K$-isometry, the transported metric in \eqref{eq:dK} is not the supremum-norm metric; for instance, $d_K(0,c\one)=3c$ for $0<c<1/2$.
We show that admissible constant shifts nevertheless recover the supremum-norm distance.
We first work with the invariant $2\norminf{h}+\int_0^1h(t)\,\dd t$ from \eqref{eq:R-invariant}, without assuming that the supremum norm and the integral are preserved separately.

\begin{lemma}\label{lem:min-invariant}
For every $h\in K_{x_0}$, we have
\[
\min_I(Th)=\min_I h.
\]
\end{lemma}
\begin{proof}
The function $h-(\min_I h)\one$ belongs to $K_{x_0}$: it is nonnegative and subtracting a nonnegative constant from $h$ can only decrease each of the two defining norm-integral expressions.
Applying Proposition~\ref{prop:constant-shifts} gives
\[
Th=T(h-(\min_I h)\one)+(\min_I h)\one.
\]
Since $T(h-(\min_I h)\one)\geq0$, this implies $\min_I(Th)\geq \min_I h$.
The same reasoning for $T^{-1}$ gives the reverse inequality.
\end{proof}

For nearby derivatives, a positive shift makes the function-value term in \eqref{eq:dK} attain its maximum at zero.
After subtracting the shift and the difference of the known values of $2\norminf{h}+\int_0^1h(t)\,\dd t$, the following quantity remains.

For $h,k\in K_{x_0}$, define
\[
Q(h,k):=\max_I(h-k)-\norminf{h}+\norminf{k}.
\]

\begin{lemma}\label{lem:Q-invariant}
For all $h,k\in K_{x_0}$, we have $Q(Th,Tk)=Q(h,k)$.
\end{lemma}
\begin{proof}
First suppose $c>0$, $h+c\one\in K_{x_0}$, and $c\geq3\norminf{h-k}$.
Since $h-k+c\one\geq0$ and $\norminf{h}-\norminf{k}+c\geq0$, the function $(\norminf{h}-\norminf{k}+c)\one+J(h-k+c\one)$ is nonnegative and nonincreasing.
Its supremum norm equals its value at zero.
Formula~\eqref{eq:dK} therefore gives
\begin{align}
d_K(h+c\one,k)&=3c+\norminf{h}-\norminf{k}+\int_0^1(h-k)(t)\,\dd t+\max_I(h-k)\notag\\
&=3c+2(\norminf{h}-\norminf{k})+\int_0^1(h-k)(t)\,\dd t+Q(h,k).\label{eq:shifted-distance}
\end{align}
This computation itself only requires $c\geq\norminf{h-k}$ and admissibility of the shift of $h$.
By \eqref{eq:metric-comparison},
\[
\norminf{Th-Tk}\leq d_K(Th,Tk)=d_K(h,k)\leq3\norminf{h-k}\leq c.
\]
Moreover, $Th+c\one\in K_{x_0}$ by Proposition~\ref{prop:constant-shifts}.
Applying \eqref{eq:shifted-distance} to the images and using \eqref{eq:R-invariant} together with $d_K(Th+c\one,Tk)=d_K(h+c\one,k)$ proves
\begin{equation}\label{eq:Q-local}
Q(Th,Tk)=Q(h,k)
\end{equation}
for such pairs.

We now pass from this local equality to all pairs.
Notice that for $h,k,u\in K_{x_0}$,
\begin{equation}\label{eq:Q-subadditive}
Q(h,k)\leq Q(h,u)+Q(u,k),
\end{equation}
by subadditivity of the maximum functional and cancellation of the intermediate supremum-norm terms.
Suppose first that $\norminf{h}+\int_{x_0}^1h(t)\,\dd t<1$ and $\norminf{k}+\int_{x_0}^1k(t)\,\dd t<1$.
Put
\[
\begin{aligned}
\eta&=1-\max\{\norminf{h}+\int_{x_0}^1h(t)\,\dd t,\norminf{k}+\int_{x_0}^1k(t)\,\dd t\}>0,\\
\rho&=2-\max\{2\norminf{h}+\int_0^1h(t)\,\dd t,2\norminf{k}+\int_0^1k(t)\,\dd t\}>0.
\end{aligned}
\]
For $h_s=(1-s)h+sk$, $0\leq s\leq1$, convexity gives $\norminf{h_s}+\int_{x_0}^1h_s(t)\,\dd t\leq1-\eta$ and $2\norminf{h_s}+\int_0^1h_s(t)\,\dd t\leq2-\rho$.
Choose
\[
c=\min\left\{\frac{\eta}{2-x_0},\frac{\rho}{6}\right\}>0.
\]
For a nonnegative function $u$ and $c\geq0$,
\[
\begin{aligned}
\norminf{u+c\one}+\int_{x_0}^1(u(t)+c)\,\dd t&=\norminf{u}+\int_{x_0}^1u(t)\,\dd t+(2-x_0)c,\\
2\norminf{u+c\one}+\int_0^1(u(t)+c)\,\dd t&=2\norminf{u}+\int_0^1u(t)\,\dd t+3c.
\end{aligned}
\]
It follows that $h_s+c\one\in K_{x_0}$ for every $s\in[0,1]$.
In particular, the same positive shift is admissible along the entire segment.
Choose $n\geq1$ such that $3\norminf{h-k}/n\leq c$ and let $h_i=(1-i/n)h+(i/n)k$ for $0\leq i\leq n$.
Equation~\eqref{eq:Q-local} applies to $h_i,h_{i+1}$.
By \eqref{eq:Q-subadditive},
\begin{align*}
Q(Th,Tk)&\leq\sum_{i=0}^{n-1}Q(Th_i,Th_{i+1})=\sum_{i=0}^{n-1}Q(h_i,h_{i+1})=Q(h,k).
\end{align*}
The final equality holds because $h_i-h_{i+1}=(h-k)/n$ and the supremum-norm terms telescope.

For arbitrary $h,k\in K_{x_0}$, set $h_j=(1-1/j)h$ and $k_j=(1-1/j)k$ for $j\geq2$.
Both satisfy the strict constraints used above.
They converge to $h,k$ in the supremum norm, and $T$ and $Q$ are continuous in that norm by \eqref{eq:metric-comparison} and the continuity of the maximum functional.
Passing to the limit gives $Q(Th,Tk)\leq Q(h,k)$ for all $h,k\in K_{x_0}$.
The argument applies equally to $T^{-1}$, since it fixes zero, preserves $2\norminf{h}+\int_0^1h(t)\,\dd t$ and satisfies Proposition~\ref{prop:constant-shifts}.
Applying the inequality for $T^{-1}$ at $Th,Tk$ proves equality.
\end{proof}

\begin{proposition}\label{prop:extrema-recovery}
For all $h,k\in K_{x_0}$,
\begin{equation}\label{eq:T-invariants}
\norminf{Th}=\norminf{h},\qquad \int_0^1(Th)(t)\,\dd t=\int_0^1h(t)\,\dd t,
\end{equation}
and
\begin{equation}\label{eq:extrema}
\max_I(Th-Tk)=\max_I(h-k),\qquad
\min_I(Th-Tk)=\min_I(h-k).
\end{equation}
In particular, $T$ is a surjective supremum-norm isometry of $K_{x_0}$.
\end{proposition}
\begin{proof}
Since $Q(0,h)=\norminf{h}-\min_I h$, Lemmas~\ref{lem:min-invariant} and \ref{lem:Q-invariant}, together with $T0=0$, yield $\norminf{Th}=\norminf{h}$.
As $2\norminf{h}+\int_0^1h(t)\,\dd t$ is preserved by \eqref{eq:R-invariant}, so is $\int_0^1h(t)\,\dd t$.
Substituting preservation of the supremum norm into the equality for $Q$ gives the maximum identity in \eqref{eq:extrema}.
Interchanging $h$ and $k$ and negating yields the minimum identity.
These two identities give $\norminf{Th-Tk}=\norminf{h-k}$.
\end{proof}

\section{Metric rigidity and linear extension}
Since $\Sp_{x_0}$ is the disjoint union of $U_{x_0}$ and $E_{x_0}$, we first establish the rigidity of each of these two sets under surjective self-isometries of $\Sp_{x_0}$.
\begin{lemma}\label{lem:right-rigidity}
Every surjective isometry $\Delta:\Sp_{x_0}\to \Sp_{x_0}$ fixes $U_{x_0}$ pointwise.
\end{lemma}
\begin{proof}
By Proposition~\ref{prop:extrema-recovery} and the density statement in Proposition~\ref{prop:right-derivative-bijection}, $T$ extends uniquely to a supremum-norm isometry $\overline T$ of $\overline{K_{x_0}}$ onto itself.
Specifically, extending both $T$ and $T^{-1}$ by continuity gives mutually inverse isometries on the closed set $\overline{K_{x_0}}$.
Since $\overline{K_{x_0}}$ is a closed convex body with nonempty interior, Mankiewicz's theorem \cite{Mankiewicz} extends $\overline T$ to a surjective affine isometry $L:\Ci\to\Ci$.
As $L0=T0=0$, $L$ is $\R$-linear.
Applying Proposition~\ref{prop:constant-shifts} to $h=0$ and $c=1/4$ gives $L(\one/4)=\one/4$, hence $L\one=\one$.

By the real Banach--Stone theorem, there exist a homeomorphism $\sigma:I\to I$ and a continuous function $w:I\to\{-1,1\}$ such that $Lh=w(h\circ\sigma)$.
Taking $h=\one$ yields $w=\one$, so
\begin{equation}\label{eq:BS}
Lh=h\circ\sigma\qquad(h\in\Ci).
\end{equation}
By \eqref{eq:T-invariants}, $L$ preserves the integral $\int_0^1h(t)\,\dd t$ on $K_{x_0}$.
Every nonnegative continuous function has a positive scalar multiple in $K_{x_0}$ (choose the scalar sufficiently small), and every continuous real function is a difference of two nonnegative ones.
By linearity, it follows that
\begin{equation}\label{eq:integral-preserved}
\int_0^1 h(\sigma(t))\,\dd t=\int_0^1 h(t)\,\dd t\qquad(h\in\Ci).
\end{equation}
For $0<s<1$, approximate the indicator of $[0,s]$ by the continuous functions $u_n(t)=\max\{0,1-n\max(t-s,0)\}$, which decrease pointwise to that indicator.
Applying \eqref{eq:integral-preserved} to $u_n$ and using dominated convergence yields $|\sigma^{-1}([0,s])|=s$, where $|\cdot|$ denotes Lebesgue measure.
An interval homeomorphism is strictly monotone.
In the increasing case the preimage has length $\sigma^{-1}(s)$; in the decreasing case it has length $1-\sigma^{-1}(s)$.
Continuity at the endpoints therefore shows that either $\sigma(x)=x$ for all $x\in I$ or $\sigma(x)=1-x$ for all $x\in I$.

To exclude reflection, take $f(x)=x/2$ and $g(x)=x^2/3$.
Both belong to $U_{x_0}$ for every $x_0\in I$.
If $\sigma(x)=1-x$, then \eqref{eq:BS} and \eqref{eq:Phi-general} give $\Delta f=f$ and $\Delta g=\widetilde g$, where
\[
\widetilde g(x)=\frac{2x}{3}-\frac{x^2}{3}.
\]
In fact, the derivative of $g$ is reflected while its terminal value remains $1-\norminf{g'}=1/3$.
The derivative differences $f'-g'$ and $f'-\widetilde g'$ both have supremum norm $1/2$.
For the function values,
\[
f-g=\frac3{16}-\frac{(x-3/4)^2}{3},
\qquad f-\widetilde g=\frac{(x-1/4)^2}{3}-\frac1{48}.
\]
On $I$ the first ranges from $0$ to $3/16$ and the second from $-1/48$ to $1/6$.
Thus
\[
\norm{f-g}=\frac{11}{16}\ne\frac23=\norm{f-\widetilde g},
\]
contradicting the isometry of $\Delta$.
We conclude that $\sigma=\id_I$, and hence $T=\id_{K_{x_0}}$.
The bijection in Proposition~\ref{prop:right-derivative-bijection} then gives $\Delta f=f$ for all $f\in U_{x_0}$.
\end{proof}

\begin{lemma}\label{lem:left-rigidity}
Every surjective isometry $\Delta:\Sp_{x_0}\to \Sp_{x_0}$ fixes $E_{x_0}$ pointwise.
\end{lemma}
\begin{proof}
Take $f\in E_{x_0}$ and put $h=f'$.
Formula~\eqref{eq:Phi-minus} gives $f(0)=\norminf{h}-1$ and $2\norminf{h}+\int_0^1h(t)\,\dd t\leq2$.
Define its right-branch counterpart by the same formula as in \eqref{eq:Phi-general}:
\[
f^+(x)=1-\norminf{h}-\int_x^1 h(t)\,\dd t.
\]
We have $(f^+)'=h$ and
\begin{equation}\label{eq:branch-difference}
f^+-f=(2-2\norminf{h}-\int_0^1h(t)\,\dd t)\one\geq0.
\end{equation}
It follows that $f^+(x_0)\geq f(x_0)\geq0$.
Since $f^+(1)=1-\norminf{h}$ and $f^+(1)\geq-f^+(0)$, we have $f^+\in \Sp_{x_0}$ and $h\in \overline{K_{x_0}}$.
If $2\norminf{h}+\int_0^1h(t)\,\dd t<2$, then $f^+\in U_{x_0}$.

For every $f\in E_{x_0}$, its counterpart $f^+$ belongs to the closure of $U_{x_0}$: the derivatives $(1-1/j)h$ belong to $K_{x_0}$ for $j\geq2$, and their images under $\Phi$ converge to $f^+$ in the sum norm by \eqref{eq:Phi-general}.
Thus $\Delta f^+=f^+$ by Lemma~\ref{lem:right-rigidity} and continuity.
Moreover,
\[
\begin{aligned}
\norm{\one-f^+}&=2\norminf{h}+\int_0^1h(t)\,\dd t,\\
\norm{f^+-f}&=2-2\norminf{h}-\int_0^1h(t)\,\dd t,\\
\norm{\one-f}&=2.
\end{aligned}
\]
Put $g=\Delta f\in E_{x_0}$.
Since $\Delta$ fixes $\one$ and $f^+$, we have $\norm{\one-g}=\norm{\one-f^+}+\norm{f^+-g}$.
Equality in the derivative component gives $\norminf{g'}=\norminf{h}+\norminf{h-g'}$.
Applying the same argument to $\Delta^{-1}$ and $g$ yields $\norminf{h}=\norminf{g'}+\norminf{g'-h}$.
Hence $g'=h$, and \eqref{eq:Phi-minus} determines the same initial value for $g$ and $f$.
This proves $g=f$, including the case $2\norminf{h}+\int_0^1h(t)\,\dd t=2$.
\end{proof}

\begin{proof}[\textbf{Proof of Theorem~\ref{thm:positive-rigidity}}]
By the proof of Lemma~\ref{lem:constant-fixed}, the diameter of $\Sp_0$ is not attained, whereas that of $\Sp_{x_0}$ is attained whenever $x_0>0$.
Hence either $x_1=x_2=0$ or both $x_1$ and $x_2$ are positive.
The metric characterizations in Lemmas~\ref{lem:constant-fixed} and \ref{lem:right-branch} yield
\[
\Delta\one=\one,\qquad \Delta(U_{x_1})=U_{x_2}.
\]
Consequently, $T=D\circ\Delta\circ\Phi:K_{x_1}\to K_{x_2}$ is a bijective $d_K$-isometry fixing zero, since the formula for $d_K$ in \eqref{eq:dK} is independent of the anchor.

Applying Lemma~\ref{lem:unique-midpoint} to the two positive spheres, the proof of Proposition~\ref{prop:constant-shifts} yields the same shift equivariance for $T$ and $T^{-1}$.
The arguments in Section~4 then remain valid, with admissible shifts and line segments taken in $K_{x_1}$ for $T$ and in $K_{x_2}$ for $T^{-1}$.
In particular, $T$ preserves the supremum-norm distance and the integral.

By the same argument as in the proof of Lemma~\ref{lem:right-rigidity}, $T$ extends to a surjective $\R$-linear isometry $L:\Ci\to\Ci$ of the form $Lh=h\circ\sigma$, where $\sigma$ is either the identity or the reflection of $I$.
Since the functions $f(x)=x/2$ and $g(x)=x^2/3$ used to exclude reflection belong to $U_{x_0}$ for every $x_0\in I$, the same contradiction excludes reflection here.
Hence $L=\id_{\Ci}$, and the surjectivity of $T$ gives $K_{x_1}=K_{x_2}$.

Suppose that $x_1<x_2$.
Choose a nonzero nonnegative continuous function $h$ supported in $(x_1,x_2)$ and normalize it so that
\[
\norminf{h}+\frac34\int_0^1h(t)\,\dd t=1.
\]
Then
\[
\norminf{h}+\int_{x_2}^1h(t)\,\dd t
=\norminf{h}<1\qquad\text{and}\qquad
2\norminf{h}+\int_0^1h(t)\,\dd t<2,
\]
so $h\in K_{x_2}$.
On the other hand,
\[
\norminf{h}+\int_{x_1}^1h(t)\,\dd t
=\norminf{h}+\int_0^1h(t)\,\dd t>1,
\]
so $h\notin K_{x_1}$, a contradiction.
The case $x_2<x_1$ is analogous.
Therefore $x_1=x_2=:x_0$.

Lemmas~\ref{lem:right-rigidity} and \ref{lem:left-rigidity} now yield $\Delta=\id_{\Sp_{x_0}}$.
The identity operator on $C^1(I)$ is an $\R$-linear extension of $\Delta$.
Any $\R$-linear extension of $\Delta$ fixes $\K_{x_0}$ pointwise, since every nonzero element of $\K_{x_0}$ is a positive multiple of an element of $\Sp_{x_0}$.
As $\K_{x_0}$ is generating, linearity gives the uniqueness of this extension.
\end{proof}

\begin{remark}
Since $f(0)\geq0$ is equivalent to $f\geq0$ whenever $f'\geq0$, we have
\[
\K_0=\{f\in C^1(I):f\geq0,\ f'\geq0\}.
\]
Thus Theorem~\ref{thm:positive-rigidity} includes the case of the positive cone consisting of all nonnegative, nondecreasing continuously differentiable functions.
\end{remark}

\begin{corollary}\label{cor:ambient-order}
Let $x_1,x_2\in I$.
If $\widetilde{\Delta}:C^1(I)\to C^1(I)$ is a surjective $\R$-linear isometry satisfying $\widetilde{\Delta}(\K_{x_1})=\K_{x_2}$, then $x_1=x_2$ and $\widetilde{\Delta}=\id_{C^1(I)}$.
\end{corollary}
\begin{proof}
The restriction of $\widetilde{\Delta}$ to $\Sp_{x_1}$ is a surjective isometry onto $\Sp_{x_2}$, so the conclusion follows from Theorem~\ref{thm:positive-rigidity}.
\end{proof}

\section*{Acknowledgments}
The author acknowledges the use of OpenAI's ChatGPT in the development of mathematical arguments and the preparation and revision of the \LaTeX{} manuscript, including the arguments for arbitrary anchors.
The basic metric invariants and the derivative-coordinate parametrization for the endpoint anchor $x_0=0$ were developed independently by the author in the original manuscript.
The author has independently verified all AI-generated mathematical arguments and bibliographic references and takes full responsibility for the correctness and content of this manuscript.

\end{document}